\documentclass[12pt]{amsart}
\usepackage{verbatim}
\usepackage{amssymb,amsmath}

\newtheorem{lemma}{Lemma} [section]
\newtheorem{thm}[lemma]{Theorem}
\newtheorem{cor}[lemma]{Corollary}

\theoremstyle{remark}
\newtheorem*{remark}{Remark}

\newcommand{\rmv}[1]{}
\newcommand{\F}{{\mathbb F}}

\DeclareMathOperator{\Gal}{Gal}
\newcommand{\Z}{{\mathbb Z}}
\newcommand{\Q}{{\mathbb Q}}

\newcommand{\C}{{\mathbb C}}

\newcommand{\tth}{^{\operatorname{th}}}


\begin{document}



\title[On some permutation polynomials]{On some permutation polynomials
over $\F_q$ of the form $x^r h(x^{(q-1)/d})$}

\author{Michael E. Zieve}
\address{Center for Communications Research, 805 Bunn Drive, Princeton NJ 08540}

\email{zieve@math.rutgers.edu}
\urladdr{http://www.math.rutgers.edu/$\sim$zieve}

\begin{abstract}
In a recent paper, Akbary and Wang gave a sufficient condition
for $x^u+x^r$ to permute $\F_q$, in terms of the period of a certain
sequence involving sums of cosines.  As an application they gave
necessary and sufficient conditions in case $u,r,q$ satisfy certain
special properties.  We show that the Akbary-Wang sufficient condition
follows from a more general sufficient condition which does not
involve sums of cosines.  This leads to vastly simpler proofs of the
Akbary-Wang results, as well as generalizations to polynomials of
the form $x^r h(x^{(q-1)/d})$.
\end{abstract}

\maketitle


\section{Introduction}

A polynomial over a finite field is called a \emph{permutation polynomial}
if it permutes the elements of the field.  These polynomials first arose
in work of Betti~\cite{B}, Mathieu~\cite{M} and Hermite~\cite{H} as a way
to represent permutations.
A general theory was developed by Hermite~\cite{H} and Dickson~\cite{D}, with
many subsequent developments by Carlitz and others.  The study of permutation
polynomials has intensified in the past few decades, due both to various
applications (e.g., \cite{golomb,powerline,DD,turbo}) and to an increasing
appreciation of the depth of the subtleties inherent to
permutation polynomials themselves (for instance, work on permutation
polynomials led to a bound on the automorphism group of a curve with
ordinary Jacobian \cite{GZ}).

The interesting aspect of permutation polynomials is the interplay
between two different ways of representing an object: combinatorially, as
a mapping permuting a set, and algebraically, as a polynomial.
This is exemplified by one of the first results in the subject, namely
that there is no permutation polynomial over $\F_q$ of degree $q-1$
if $q>2$ \cite{H}.  Much recent work has focused on low-degree
permutation polynomials, as these have quite remarkable properties:
for instance, a polynomial of degree at most $q^{1/4}$ which permutes $\F_q$
will automatically permute $\F_{q^n}$ for infinitely many~$n$.
The combined efforts of several mathematicians have led to a handful of
families of such polynomials, and to an avenue towards proving there are no
others \cite{D,DL,Co70,FGS,Muller,CM,LZ,GM,GMZ,GRZ}.

A different line of research focuses not on the degree of a permutation
polynomial but instead on the number of terms.  The simplest class of
nonconstant polynomials are the monomials $x^n$ with $n>0$, and one easily
checks that $x^n$ permutes
$\F_q$ if and only if $n$ is coprime to $q-1$.  However, for
binomials the situation becomes much more mysterious.  Despite the attention
of numerous authors since the 1850's (cf., e.g., \cite{B,M,H,Br,C,NR,T,S,W2,KL,
Wang,AW,MZ2}),
the known results seem far from telling
the full story of permutation binomials.  This brings us to the present paper.
In the recent paper \cite{AW},
Akbary and Wang considered binomials of the form $f(x) = x^u+x^r$
with $u>r>0$.  They gave sufficient conditions for $f$ to permute
$\F_q$ in terms of the period of the sequence ($a_n$~mod~$p$),
where $p$ is the characteristic of $\F_q$ and, with
$d:=(q-1)/\gcd(q-1,u-r)$,
\[
a_n := \sum_{t=1}^{\frac{d-1}2} \left(2\cos\frac{\pi (2t-1)}
d\right)^n.
\]
(One can show that every $a_n$ is an integer.)

As an application, they gave necessary and sufficient conditions
for $x^u+x^r$ to permute $\F_q$ in the two special cases
\begin{enumerate}
\item $p\equiv 1\pmod{d}\,\text{ and }\, d\mid \log_p q$.
\item $p\equiv -1\pmod{d}$.
\end{enumerate}

The proofs in \cite{AW} relied on facts about the coefficients of Chebychev
polynomials, Hermite's criterion, properties of recursive sequences, 
lacunary sums of binomial coefficients, and various unpublished results
about factorizations of Chebychev polynomials, among other things.
In this paper we give vastly shorter and simpler proofs which avoid all these
ingredients, and which yield more general results.
In particular, we will show that the sequence $a_n$ does not play an
essential role for these results, and in fact stating results in terms of
$a_n$ obscures the essence of the situation.

We will prove the following sufficient condition for permutation binomials,
in which (for any $d>0$) $\mu_d$ denotes the set of $d^{\operatorname{th}}$
roots of unity in the algebraic closure of $\F_q$:

\begin{thm}
\label{gen}
Pick $u>r>0$ and $a\in\F_q^*$.  Write $s:=\gcd(u-r,q-1)$ and
$d:=(q-1)/s$.
Suppose that $(\eta+a/\eta)\in\mu_s$ for every $\eta\in\mu_{2d}$.
Then $x^u+ax^r$ permutes\/ $\F_q$ if and only if
$-a\notin\mu_d$ and $\gcd(r,s)=1$ and $\gcd(2d,u+r)\le 2$.
\end{thm}

We emphasize that this condition applies to arbitrary binomials, unlike
the condition in \cite{AW} which only applied to binomials with both
coefficients being~$1$.  Superficially the condition in \cite{AW}
looks quite different from Theorem~\ref{gen}, since the former requires a
constraint on the period of $(a_n \text{ mod }p)$; however, in
Section~\ref{sec aw} we will show that
the hypotheses of Theorem~\ref{gen} are satisfied whenever the hypotheses
of \cite[Thm.~1.1]{AW} are satisfied.  We note further that, for
both theoretical and practical purposes, our hypotheses are easier to test
than those in \cite{AW}.

In the forthcoming paper \cite{AAW}, the two families of permutation binomials
from \cite{AW} are generalized to families of permutation polynomials of the form
$x^r(1+x^s+x^{2s}+\dots+x^{ks})$, with similar proofs to those in \cite{AW}.
We now exhibit two vastly more general families of permutation polynomials which
include the polynomials from \cite{AW} and \cite{AAW} as quite special cases.

\begin{thm}
\label{introlaigle}
Let $d,r>0$ satisfy $d\mid (q-1)$.
Suppose that $q=q_0^m$ where $q_0\equiv 1\pmod d$ and $d\mid m$,
and pick $h\in\F_{q_0}[x]$.  
Then $f(x):=x^r h(x^{(q-1)/d})$
permutes $\F_q$ if and only if $\gcd(r,(q-1)/d)=1$ and $h$ has no
roots in $\mu_d$.
\end{thm}

This is equivalent to a forthcoming result of Laigle-Chapuy \cite{LC}; our
proof is significantly simpler than
that in \cite{LC}.  The first class of permutation binomials
from \cite{AW} is the special case that $q_0=p$ and $h=x^e+1$, where $\gcd(e,d)=1$.

In our next result we use the notation $h_k(x):=x^{k-1}+x^{k-2}+\dots+1$.

\begin{thm}
\label{introneg}
Pick integers $t\ge 0$ and $r,v,k,\ell>0$, and put $s:=\gcd(q-1,v)$,
$d:=(q-1)/s$, and $d_0:=d/\gcd(d,\ell-1)$.  Suppose that $q=q_0^m$,
where $m$ is even and $q_0\equiv -1\pmod{d}$.
Pick $\hat h\in\F_{q_0}[x]$ and let $h:= h_k(x)^t \hat h(h_{\ell}(x)^{d_0})$.
Then $f:=x^r h(x^v)$ permutes $\F_q$ if and only if
$\gcd(r,s)=1$,\, $\gcd(2r+(k-1)tv,2d)=2$ and $h$ has no roots in $\mu_d$.
\end{thm}

The second class of permutation binomials from \cite{AW} is the special case
that $q_0=p$ and $h=h_2$.

{\bf Notation:} Throughout this paper, $q$ is a power of the prime $p$,
and $\mu_d$ denotes the set of $d\tth$ roots of unity in the algebraic
closure of $\F_q$.  Also, $h_k(x):=x^{k-1}+x^{k-2}+\dots+1$.

\section{Proofs}

We begin with a simple lemma reducing the question whether a polynomial
permutes $\F_q$ to the question whether a related polynomial permutes
a particular subgroup of $\F_q^*$.

\begin{lemma}
\label{lwl}
Pick $d,r>0$ with $d\mid (q-1)$, and let $h\in\F_q[x]$.
Then $f(x):=x^r h(x^{(q-1)/d})$ permutes\/ $\F_q$ if and only if
both
\begin{enumerate}
\item $\gcd(r,(q-1)/d)=1$\, and
\item $x^r h(x)^{(q-1)/d}$ permutes $\mu_d$.
\end{enumerate}
\end{lemma}

\begin{proof}
Write $s:=(q-1)/d$.
For $\zeta\in\mu_s$, we have $f(\zeta x) = \zeta^r f(x)$.
Thus, if $f$ permutes $\F_q$ then $\gcd(r,s)=1$.
Conversely, if $\gcd(r,s)=1$ then the values of $f$ on $\F_q$
consist of all the $s^{\operatorname{th}}$ roots of the values
of
\[
f(x)^s = x^{rs} h(x^s)^s.
\]
But the values of $f(x)^s$ on $\F_q$ consist of $f(0)^s=0$
and the values of $g(x):=x^r h(x)^s$ on $(\F_q^*)^s$.
Thus, $f$ permutes $\F_q$ if and only if $g$ is bijective
on $(\F_q^*)^s=\mu_d$.
\end{proof}

\begin{remark}
A more complicated criterion for $f$ to permute $\F_q$ was given
by Wan and Lidl~\cite[Thm.\ 1.2]{WL}.
\end{remark}

The difficulty in applying Lemma~\ref{lwl} is verifying condition (2).
Here is one situation where this is easy:

\begin{cor}
\label{apply}
Pick $d,r,n>0$ with $d\mid (q-1)$, and let $h\in\F_q[x]$.
Suppose $h(\zeta)^{(q-1)/d}=\zeta^n$ for every $\zeta\in\mu_d$.  Then
$f(x):=x^r h(x^{(q-1)/d})$ permutes\/ $\F_q$ if and only if
$\gcd(r+n,d)=\gcd(r,(q-1)/d)=1$.
\end{cor}

Our next results give choices for the parameters satisfying the
hypotheses of Corollary~\ref{apply}.

\begin{thm}
\label{laigle}
Let $d,r>0$ satisfy $d\mid (q-1)$.
Suppose that $q=q_0^m$ where $q_0\equiv 1\pmod d$ and $d\mid m$,
and pick $h\in\F_{q_0}[x]$.  
Then $f(x):=x^r h(x^{(q-1)/d})$
permutes $\F_q$ if and only if $h$ has no roots in $\mu_d$
and $\gcd(r,(q-1)/d)=1$.
\end{thm}

\begin{proof}
We may assume $\gcd(r,(q-1)/d)=1$, since otherwise $f$ does not permute
$\F_q$ (by Lemma~\ref{lwl}).
Since $q_0\equiv 1\pmod d$, we have
\[
\frac{q_0^d-1}{q_0-1} = \sum_{i=0}^{d-1}q_0^i\equiv 0\pmod{d}.
\]
Hence $q_0-1$ divides $(q_0^d-1)/d$, which divides $(q-1)/d$;
since $d\mid(q_0-1)$, it follows that $d$ divides $(q-1)/d$, so
since $\gcd(r,(q-1)/d)=1$ we have $\gcd(r,q-1)=1$.

For $\zeta\in\mu_d$ we have $\zeta\in\F_{q_0}$, so
$h(\zeta)\in\F_{q_0}$.  Since $f(0)=0$, if $f$ permutes $\F_q$
then $h(\zeta)\ne 0$.  Conversely, if $h(\zeta)\ne 0$ then
(since $q_0-1$ divides $(q-1)/d$) we have
$h(\zeta)^{(q-1)/d}=1$.  Now the result follows from
Corollary~\ref{apply} (with $n=d$).
\end{proof}

\begin{remark}
Theorem~\ref{laigle} is a reformulation of a result from \cite{LC},
which contains a different proof.  (Note that \cite[Thm.~4.3]{LC} is
false, a counterexample being $P=x^3+x$ over $\F_3$; to correct it one
should remove the polynomials $P$.)
\end{remark}

We now exhibit some polynomials $h$ for which we can determine when $h$ has
roots in $\mu_d$.

\begin{cor}
\label{special}
Pick positive integers $d,e,r,k,t$ with $d\mid(q-1)$ and $\gcd(d,e)=1$.
Suppose that $q=q_0^m$ where $q_0\equiv 1\pmod{d}$ and $d\mid m$.
Then $f(x):=x^r h_k(x^{e(q-1)/d})^t$ permutes $\F_q$ if and only if
$\gcd(k,pd)=\gcd(r,(q-1)/d)=1$.
\end{cor}

\begin{remark}
The case that $q_0=p$, $k=2$, and $t=1$ was treated in \cite{AW}.
The case that $q_0=p$, $t=e=1$, and both $q$ and $d$ are odd was
treated in \cite{AAW}.
The results in both \cite{AW} and \cite{AAW} involved the superfluous
condition $\gcd(2r+(k-1)es,d)=1$.
\end{remark}

\begin{thm}
\label{neg}
Pick integers $t\ge 0$ and $r,v,k,\ell>0$, and put $s:=\gcd(q-1,v)$, $e:=v/s$,
$d:=(q-1)/s$, and $d_0:=d/\gcd(d,\ell-1)$.  Suppose that $q=q_0^m$,
where $m$ is even and $q_0\equiv -1\pmod{d}$.
Pick $\hat h\in\F_{q_0}[x]$, and let
$h(x):= h_k(x)^t \hat h(h_{\ell}(x)^{d_0})$.
Then $f:=x^r h(x^v)$ permutes $\F_q$ if and only if
$\gcd(r,s)=1$,\, $\gcd(2r+(k-1)tv,2d)=2$ and $h$ has no roots in $\mu_d$.
\end{thm}

\begin{proof}
Our hypotheses imply the divisibility relations
\[
q_0-1 = \frac{q_0^2-1}{q_0+1} \left|\, \frac{q-1}{q_0+1} \,\right| \frac{q-1}d
=s.
\]
We may assume $h(x^e)$ has no roots in $\mu_d$, since otherwise
Lemma~\ref{lwl} implies $f$ does not permute $\F_q$; since $\gcd(d,e)=1$,
this says $h$ has no roots in $\mu_d$.  Hence $\hat h(h_{\ell}(x)^{d_0})$ has
no roots in $\mu_d$.
For $\zeta\in\mu_d\setminus\mu_1$, the hypothesis $d\mid(q_0+1)$ implies
$\zeta^{q_0}=1/\zeta$, so
\[
h_{\ell}(\zeta)^{q_0} = \left(\frac{\zeta^{\ell}-1}{\zeta-1}\right)^{q_0}
= \frac{\zeta^{-\ell}-1}{\zeta^{-1}-1} 
= \frac{h_{\ell}(\zeta)}{\zeta^{\ell-1}};
\]
hence $h_{\ell}(\zeta)^{d_0q_0} = h_{\ell}(\zeta)^{d_0}$, so
$h_{\ell}(\zeta)^{d_0}\in\F_{q_0}$.  Also $h_{\ell}(1)\in\F_{q_0}$.
Thus, for any $\zeta\in\mu_d$ we have
$\hat h(h_{\ell}(\zeta^e)^{d_0})\in\F_{q_0}^*$.
Since $(q_0-1)\mid s$, we conclude that $h(\zeta^e)^s=h_k(\zeta^e)^{ts}$.
As above, $h_k(\zeta)^{t(q_0-1)}=1/\zeta^{t(k-1)}$, so
$h(\zeta^e)^s = 1/\zeta^{e(k-1)ts/(q_0-1)}$, whence the
result follows from Corollary~\ref{apply}.
\end{proof}

\begin{remark}
There would be counterexamples to Theorem~\ref{neg} if we did not
require $m$ even; such examples necessarily have $d=2$.  Also,
Theorem~\ref{neg} immediately generalizes to the case that $h$ is the
product of several polynomials of the same shapes as the two factors of
$h$ described in the theorem, and moreover
we may replace $h$ by any polynomial congruent to it modulo $x^d-1$.
\end{remark}

\begin{cor}
\label{specialneg}
Pick positive integers $t,d,e,r,k$ with $d\mid(q-1)$ and $\gcd(d,e)=1$,
and put $s:=(q-1)/d$.
Suppose that $q=q_0^m$ where $m$ is even and $q_0\equiv -1\pmod{d}$.
Then $f(x):=x^r h_k(x^{e(q-1)/d})^t$ permutes $\F_q$ if and only if
$\gcd(r,s)=\gcd(k,pd)=1$ and $\gcd(2r+(k-1)tes,2d)=2$.
\end{cor}

\begin{remark}
The hypotheses of Corollary~\ref{specialneg} are satisfied whenever $d$ is an
odd prime divisor of $q-1$ such that $p$ has even order modulo~$d$.
The case that $d=7$, $t=1$, and $k=2$ was treated in \cite{AW0}, although the
result in \cite{AW0} includes the superfluous condition $2^s\equiv 1\pmod{p}$.
The case that $q_0=p$, $t=1$, and $k=2$ was treated in \cite{AW}.
The case that $q_0=p$, $t=e=1$, and both $q$ and $d$ are odd was treated
in \cite{AAW}.
\end{remark}

Now we prove a general sufficient criterion for permutation
binomials:
\begin{thm}
\label{bin}
Pick $u>r>0$ and $a\in\F_q^*$.  Write $s:=\gcd(u-r,q-1)$ and
$d:=(q-1)/s$.
Suppose that $(\eta+a/\eta)\in\mu_s$ for every $\eta\in\mu_{2d}$.
Then $x^u+ax^r$ permutes\/ $\F_q$ if and only if
$-a\notin\mu_d$ and $\gcd(r,s)=1$ and $\gcd(2d,u+r)\le 2$.
\end{thm}

\begin{proof}
Write $e:=(u-r)/s$, so $\gcd(e,d)=1$.
By Lemma~\ref{lwl}, $f(x):=x^u+ax^r$ permutes $\F_q$ if and only if
$\gcd(r,s)=1$ and $g(x):=x^r (x^e+a)^s$ permutes $\mu_d$.  In
particular, if $x^u+ax^r$ permutes $\F_q$ then $g$ has no roots
in $\mu_d$, or equivalently $-a\notin\mu_d$.  Henceforth we assume
$\gcd(r,s)=1$ and $-a\notin\mu_d$, so $f$ permutes $\F_q$ if and only
if $g$ is injective on $\mu_d$.  This condition is equivalent to
injectivity of $g(x^2)$ on $\mu_{2d}/\mu_2$.  But for $\eta\in\mu_{2d}$
we have
\begin{align*}
g(\eta^2) &= \eta^{2r} (\eta^{2e}+a)^s \\
&= \eta^{2r+es}\left(\eta^e + \frac{a}{\eta^e}\right)^s \\
&= \eta^{2r+es}.
\end{align*}
Finally, $x^{2r+es}$ is injective on $\mu_{2d}/\mu_2$ if and only if
$\gcd(2r+es,2d)\le 2$; since $2r+es=u+r$, this completes the proof.
\end{proof}

Theorem~\ref{bin} can be generalized (with the same proof) to
polynomials with more terms:

\begin{thm}
Pick $r,e,d,t>0$ where $d\mid (q-1)$ and $\gcd(e,d)=1$.
Put $h=x^t \hat h(x^d)$ where $\hat h\in\F_q[x]$.  Pick $a\in\F_q^*$.
Suppose that every $\eta\in\mu_{d\gcd(2,d)}$ satisfies both
$\eta+a/\eta\in\mu_{t(q-1)/d}$ and $\hat h((\eta^{2e}+a)^d)\in
\mu_{(q-1)/d}$.
Then $f(x):=x^r h(x^{e(q-1)/d}+a)$ permutes $\F_q$ if and
only if $\gcd(2r+et(q-1)/d,d)=1$ and $\gcd(r,(q-1)/d)=1$.
\end{thm}

\section{Permutation binomials and generalized Lucas sequences}
\label{sec aw}

In this section we explain how our sufficient condition for permutation
binomials (Theorem~\ref{bin}) implies the analogous condition from \cite{AW},
namely \cite[Thm.~1.1]{AW}.  Some
preliminary steps are needed in order to state \cite[Thm.~1.1]{AW}.

It is easy to show that if $f(x):=x^r+x^u$ (with $0<r<u$)
permutes $\F_q$ then $f(x) = x^r(1+x^{es})$, where
\begin{equation}
\tag{$*$}
sd= q-1,\, \gcd(r,s)=\gcd(e,d)=1,\, d\text{ odd},
\text{ and } r,e,s>0.
\end{equation}
Conversely, with $p$ denoting the characteristic of $\F_q$,
\cite[Thm.~1.1]{AW} says
\begin{thm}
\label{AWthm}
For $q,s,d,r,e$ as in \emph{($*$)}, the binomial $f(x)=x^r(1+x^{es})$
permutes\/ $\F_q$ if $\gcd(2r+es,d)=1$, $2^s\equiv 1\pmod{p}$, and
the sequence
\[
a_n := \sum_{t=1}^{\frac{d-1}2} \left(2\cos\frac{\pi (2t-1)}
d\right)^n
\]
consists of integers satisfying
$a_n\equiv a_{n+s}\pmod{p}$ for every $n\ge 0$.
\end{thm}

Suppose the hypotheses of this result are satisfied, and put
$\zeta=\exp(\pi i/d)$.  Then
\[
2a_n = 2\sum_{t=1}^{\frac{d-1}2} \left(\zeta^{2t-1} + \frac{1}
{\zeta^{2t-1}}\right)^n = \sum_{\substack{\eta\in\C\setminus\{-1\} \\
\eta^d=-1}} \left(\eta + \frac{1}{\eta}\right)^n.
\]
Note that the hypotheses of Theorem~\ref{AWthm} imply $q$ odd
(since $s>0$ and $2^s\equiv 1\pmod{p}$).  Also, we now see that
$a_n\in\Z[\zeta]$ and that $a_n$ is fixed by every element of
$\Gal(\Q(\zeta)/\Q)$, so $a_n\in\Q\cap\Z[\zeta]=\Z$.
Let $\hat\zeta$ denote a fixed primitive $(2d)\tth$ root of unity
in $\F_q$, and let $\psi$ be the homomorphism $\Z[\zeta]\mapsto\F_q$
which maps $\zeta\mapsto\hat\zeta$.  Then $\psi(a_n)\equiv a_n\pmod{p}$,
so the condition $a_n\equiv a_{n+s}\pmod{p}$ is equivalent to
\[
\sum_{\substack{\eta\in\F_q\setminus\{-1\}\\ \eta^d=-1}} 
\left(\left(\eta + \frac{1}{\eta}\right)^s - 1 \right)\cdot
\left(\eta + \frac{1}{\eta}\right)^n  = 0.
\]
This condition holds for all $n\ge 0$ if and only if
\[
\sum_{\substack{\eta\in\F_q\setminus\{-1\} \\ \eta^d=-1}}
\left(\left(\eta + \frac{1}{\eta}\right)^s - 1 \right)\cdot
P\Bigl(\eta+\frac{1}{\eta}\Bigr) = 0
\]
for every $P\in\F_q[x]$.
Pick representatives $\eta_1,\eta_2,\dots,\eta_{(d-1)/2}$
for the equivalence classes of $\mu_{2d}\setminus
(\mu_d\cup\mu_2)$ under the equivalence relation
$\eta\sim 1/\eta$.  Then the values $\eta_i+1/\eta_i$ are distinct
elements of $\F_q$, so there are polynomials $P\in\F_q$
taking any prescribed values at all the $\eta_i+1/\eta_i$.
In particular, choosing $P$ to be zero at all but one of
these elements, it follows that
\begin{equation}\label{u}
\left(\eta + \frac{1}{\eta}\right)^s = 1
\end{equation}
for every $\eta$ such that $\eta^d=-1$ but $\eta\ne -1$.
The hypotheses of Theorem~\ref{AWthm} imply that $s$ is even
and $2^s\equiv 1\pmod{p}$, so (\ref{u}) holds for $\eta=-1$.
Moreover, since $d$ odd and $s$ even, the fact that (\ref{u})
holds when $\eta^d=-1$ implies that (\ref{u}) holds when
$\eta^d=1$ as well.

Thus, whenever the hypotheses of Theorem~\ref{AWthm} hold,
we will have $(\eta+1/\eta)^s=1$ for every $\eta\in\mu_{2d}$.
Since the latter is precisely the hypothesis of Theorem~\ref{bin}
in case $a=1$,
we see that Theorem~\ref{bin} implies Theorem~\ref{AWthm}.



\end{document}